\documentclass[12pt,leqno]{article}  % 建檔日期 2015年5月18日
\usepackage{amsmath,amssymb,bbm,accents}
\usepackage{amsthm,array,mathdots}
\usepackage{titlesec} % 重設章節標題所需的巨集
\usepackage[all,cmtip]{xy}
\usepackage{blkarray,arydshln} % 矩陣分隔線所需的巨集

\usepackage{tikz}

\usepackage{picins}

\usepackage{graphics}

\titlespacing{\section}{0cm}{3.5pc}{1.5pc}

\makeatletter
\def\@citex[#1]#2{\if@filesw\immediate\write\@auxout{\string\citation{#2}}\fi
  \def\@citea{}\@cite{\@for\@citeb:=#2\do
    {\@citea\def\@citea{\@citesep}\@ifundefined
       {b@\@citeb}{{\bf ?}\@warning
       {Citation `\@citeb' on page \thepage \space undefined}}%
{\csname b@\@citeb\endcsname}}}{#1}}
\def\@citesep{; }
\makeatother

\newtheoremstyle{Zhu}{}{}{\itshape}{}{\bf}{}{.5em}{}
\theoremstyle{Zhu}
\newtheoremstyle{Zhu}{}{}{\itshape}{}{\bf}{}{.5em}{}
\theoremstyle{Zhu}
\newtheorem{theorem}{Theorem}[section]

\newtheorem{lemma}[theorem]{Lemma}
\newtheorem{corollary}[theorem]{Corollary}
\newtheorem{Proposition}[theorem]{Proposition}

\newtheoremstyle{Kremark}{}{}{}{}{\bf}{}{.5em}{}
\theoremstyle{Kremark}
\newtheorem*{Remark}{Remark.}
\newtheorem{Definition}[theorem]{Definition}
\newtheorem{Example}[theorem]{Example}
\newtheorem{other}{}
\newtheorem{Note}[theorem]{Note}

\allowdisplaybreaks[1]  % 允許多行數學式

\newenvironment{pf}{\medskip\par\noindent{\bf Proof\/}.}{\hfill
$\Box$}

\title{Projective dimension and  regularity of  path  ideals  of  cycles}
\author{\begin{minipage}{0.8\textwidth}
 \hspace{4.5cm}   Guangjun Zhu \\[2mm] %\normalsize
\end{minipage} }

\vspace{-3cm}

\date{}

\begin{document}

\maketitle
\footnote{Supported by the National Natural Science Foundation of
China (11271275) and by  Foundation of Jiangsu Overseas Research \& Training Program for
University Prominent Young \& Middle-aged Teachers and Presidents and
 by Foundation of the Priority Academic Program Development of Jiangsu Higher Education Institutions. }

%\footnote{School of Mathematical Sciences,
%Soochow University,
%Suzhou, China\\
% E-mail: zhuguangjun@suda.edu.cn}

\vspace{-2cm}
\begin{abstract}
{\noindent\bf Abstract.} In this paper, we give a formula to compute all the top degree graded Betti numbers of the path ideals
of a cycle. As a consequence we can  give a formula to compute its projective dimension and regularity.

\vspace{3mm}

\vspace{3mm}
{\noindent\bf Mathematics Subject Classification ($2010$)}: 13D02; 13F55; 13C15; 13D99.

\end{abstract}

%------------------------------------S1

\section{Introduction }

\vspace{3mm}In this paper, we study some path ideals of cycle graph, these ideals can be seen as generalization
of path ideals for cycle, as studied in \cite{AF1}. Let $\Gamma=(V,E)$ be a directed
graph on vertex set $V=\{x_{1},\dots,x_{n}\}$ with edge set $E$ and let $K$ be an arbitrary
field. Consider the polynomial ring $R=K[x_{1},\dots, x_{n}]$, where we identify the
vertices of $\Gamma$ with the variables of $R$, the path ideal of  $\Gamma$ of length $m$  is the
monomial ideal
$$J_{m}(\Gamma)=
(\{x_{i_{1}}\cdots x_{i_{m}}\ |\ x_{i_{1}},\dots, x_{i_{m}}\
\mbox{is a path of length}\ m\  \mbox{in}\ \Gamma
\})\subseteq K[x_{1},\dots, x_{n}].$$

Path ideals have been  first introduced  by Conca and De Negri \cite{CN} and they generalize
arbitrary edge ideals of graphs  \cite{V1}. Since then, path ideals have attracted the
attention of a lot of researchers and they are fairly well-studied for special classes
of graphs such as the line graph and the cycle \cite{AF1,AF2,Z} and also for rooted trees \cite{BHO,CN,HT3}.
Conca and Negri in \cite{CN} showed that path ideals of directed  trees have normal and
Cohen-Macaulay Rees rings.  Restuccia and  Villarreal in \cite{RV} showed that
path ideals of complete bipartite graphs are normal, while   path
ideals of cycles are shown to have linear type in \cite{BS}. In \cite {HT3}, He and  Tuyl
study $J_{m}(\Gamma)$ in the special case that  $\Gamma$  is the line graph $L_{n}$.
They provided an exact formula for the projective dimension
of $J_{m}(L_{n})$ in terms of $m$ and $n$. They showed that:

 \begin{theorem}(Theorem 4.1)
 Let $p,m,n,d$
be  integers such that  $n=p(m+1)+d$, where $p\geq 0$, $0\leq d\leq m$ and $2\leq m\leq n$. Then
the projective dimension of $J_{m}(L_{n})$ is given by
$$pd\,(J_{m}(L_{n}))=\left\{\begin{array}{ll}
2p-1& \mbox{if}\ \ d\neq m \\
2p& \mbox{if}\ \  d=m.
\end{array}\right.$$
\end{theorem}

\vspace{3mm} Using purely combinatorial arguments, Alilooee and Faridi in \cite{AF2} also gave the above formula for
 projective dimension of $J_{m}(L_{n})$. Furthermore,  they  provide an explicit formula for
 regularity  of $J_{m}(L_{n})$ in terms of $m$ and $n$. They showed that:

 \begin{theorem}(Corollary 4.14)
 Let $p,m,n,d$
be  integers such that  $n=p(m+1)+d$, where $p\geq 0$, $0\leq d\leq m$ and $2\leq m\leq n$. Then
the regularity of $J_{m}(L_{n})$ is given by
$$reg\,(J_{m}(L_{n}))=\left\{\begin{array}{ll}
p(m-1)+1& \mbox{if}\ \ d\neq m\\
p(m-1)+m& \mbox{if}\ \ d=m.
\end{array}\right.$$
\end{theorem}

\vspace{3mm}In \cite{AF1}, using  homological methods, Alilooee and Faridi  gave exact formulas for
 projective dimension and the regularity  of $J_{m}(C_{n})$
of the cycle graph $C_{n}$. They showed that:

\begin{theorem}(Corollary 5.5)
 Let $p,m,n,d$
be  integers such that  $n=p(m+1)+d$, where $p\geq 0$, $0\leq d\leq m$ and $2\leq m\leq n$. Then

(1) The projective dimension of the path ideal of the cycle graph $C_{n}$ is given by
$$pd\,(J_{m}(C_{n}))=\left\{\begin{array}{ll}
2p&\mbox{if}\ \  d\neq 0\\
2p-1&\mbox{if}\ \  d=0.
\end{array}\right.$$

(2)  The regularity of the path ideal of the cycle graph $C_{n}$ is given by
$$reg\,(J_{m}(C_{n}))=\left\{\begin{array}{ll}
p(m-1)+d=n-2p& \mbox{if}\ \ d\neq 0\\
p(m-1)+1=n-2p+1& \mbox{if}\ \ d=0.
\end{array}\right.$$
\end{theorem}

\vspace{3mm}We generalize the notion of path ideal as the following:
Let $\Gamma$ be  a directed graph,  path ideal of $\Gamma$  of length $m$ is the  monomial ideal generated by some  paths of length $m$ in $\Gamma$, i.e.,
$I_{m}(\Gamma)=(u_{1},\dots,u_{k}),\  \mbox{where}\  u_{1},\dots,u_{k}
\ \mbox{ are some paths  of length}\  m\ \mbox{in}\ \Gamma.$
When $u_{1},\dots,u_{k}$ are all paths of length $m$ in $\Gamma$, $I_{m}(\Gamma)=J_{m}(\Gamma)$.

To the best of our knowledge,
little  is known about these ideals.  It is, therefore, of interest to determine
algebraic properties of the ideals $I_{m}(\Gamma)$. In \cite {Z},
using the notion of  Betti-splitting, Zhu study the path ideal $I_{m,l,k}=(u_{1},\dots,u_{k})$ where $u_{i}=\prod\limits_{j=1}^{m}x_{(i-1)l+j}$ for any $1\leq i\leq k$ and $l$ is an integer such that $1\leq l<m$.
She provided some exact formulas for the projective dimension and regularity
of $I_{m,l,k}$ in terms of $m$ and $n$. She showed that:

\vspace{3mm}\begin{theorem}(Theorem 3.5)
Let $k,l, m, n$ be integers such that $n=kl+(m-l)$ where  $k\geq 1$,  $m\geq 2$ and $l<\lceil \frac{m}{2}\rceil$, here $\lceil \frac{m}{2}\rceil$
denotes the smallest integer $\geq \frac{m}{2}$.
Let $I_{m,l,k}=(u_{1},\dots,u_{k})$ with $u_{i}=\prod\limits_{j=1}^{m}x_{(i-1)l+j}$ for any $1\leq i\leq k$. Then $pd\,(I_{m,l,k})=k-1$,  $reg\,(I_{m,l,k})=(k-1)(l-1)+m$.
\end{theorem}

\vspace{3mm}\begin{theorem}(Theorem 3.7)
Let $k,l, m, n$ be integers such that $n=kl+(m-l)$ where  $k\geq 1$,  $m\geq 2$ and $\lceil \frac{m}{2}\rceil\leq l<m$.
Let $I_{m,l,k}=(u_{1},\dots,u_{k})$ with $u_{i}=\prod\limits_{j=1}^{m}x_{(i-1)l+j}$ for any $1\leq i\leq k$.
If $m\equiv 0\,(\mbox{mod}\ l)$
 and we can write $n$ as $n=p(m+l)+d$ where $0\leq d<m+l$, then
 \begin{itemize}
\item[(1)]$pd\,(I_{m,l,k})=\left\{\begin{array}{ll}
2p-1& \mbox{if}\ \ d\neq m\\
2p& \mbox{if}\ \ d=m.
\end{array}\right.$
\item[(2)] $reg\,(I_{m,l,k})=\left\{\begin{array}{ll}
p(m+l-2)+1& \mbox{if}\ \ d\neq m\\
p(m+l-2)+m& \mbox{if}\ \ d=m.
\end{array}\right.$
\end{itemize}
\end{theorem}

\vspace{3mm}\begin{theorem}(Theorem 3.10)
Let $k,l, m, n$ be integers such that  $n=kl+(m-l)$ where  $k\geq 1$,  $m\geq 2$ and  $\lceil \frac{m}{2}\rceil\leq l<m$.
Let $I_{m,l,k}=(u_{1},\dots,u_{k})$ with $u_{i}=\prod\limits_{j=1}^{m}x_{(i-1)l+j}$ for any $1\leq i\leq k$.
If $m\equiv s\,(\mbox{mod}\ l)$ with $1\leq s<l$
 and we can write $n$ as $n=p(m+l-s)+d$ where $0\leq d<m+l-s$, then  $$pd\,(I_{m,l,k})=\left\{\begin{array}{ll}
2p-1& d\neq m\\
2p& d=m.
\end{array}\right.
$$
\end{theorem}

In this article we shall focus on the path ideals of an $n$-cycle graph $C_{n}$.
By promoting Alilooee and Faradi's techniques, we generalize formulas for projective dimension and regularity of the path ideal
 of the cycle graph obtained in \cite{AF1,AF2}.  we obtain that

\vspace{3mm}\begin{theorem}(Corollary 5.4)
Let $l,m,n$ be positive integers such that $2\leq m\leq n$, $l|n$ and $l\leq min\{m-1,\frac{n}{2}\}$,
let  $m\equiv s\ \mbox{mod}\ l$ with $0\leq s<l$, $t=\frac{m-s}{l}$ and $k=\frac{n}{l}$.  Suppose that
$\Delta_{m,l}(C_{n})$ $=\langle F_{1},\dots, F_{k}\rangle$ is the  path complex of the cycle  $C_{n}$ with standard labeling and $I_{m,l}(C_{n})=\mathcal{I}(\Delta_{m,l}(C_{n}))$ is the facet ideal of $\Delta_{m,l}(C_{n})$.
We can write $k$ as  $k=p(t+1)+d$, where $p\geq 0$ and $0\leq d\leq t$. Then

(1) The projective dimension of $R/I_{m,l}(C_{n})$   is given by
$$ pd\,(R/I_{m,l}(C_{n}))=\left\{\begin{array}{ll}
2p,&\mbox{if}\ d=0\\
2p+1,&\mbox{if}\ d\neq 0\\
\end{array}\right.$$

(2) The regularity of $R/I_{m,l}(C_{n})$   is given by
$$ reg\,(R/I_{m,l}(C_{n}))=\left\{\begin{array}{ll}
n-2p,&\mbox{if}\ d=0\\
n-2p-1,&\mbox{if}\ d\neq 0.\\
\end{array}\right.$$
\end{theorem}

\vspace{3mm}Our paper is organized as follows. In Section $2$, we recall some notation and
basic algebraic and combinatorial concepts used in other next chapters. In Section $3$,
we study the connected components of path ideals which will provide us with the
key to our homological computations later in Section $4$. Section $5$ is where we
apply the homological results of Section $4$ along to give a criterion to determine
all nonzero Betti numbers and projective dimension and regularity of path ideals of cycles. While
working on this paper the computer algebra systems CoCoA \cite{C} and Macaulay$2$ \cite{GS}
were used to test examples. We acknowledge the immense help that they have
provided us in this project.

\hspace{5mm}

\section{Preliminaries}

\hspace{3mm}In this section, we provide some background and fix some notation, which will
be used throughout this article. However, for more details we refer the reader to \cite{F1,F2,HH,M}.

\begin{Definition} \label{def1}
A simplicial complex $\Delta$ on vertex set   $V$ is a collection of subsets of  $V$  with the property that if $F\in \Delta$
then all  subsets of $F$ are also in $\Delta$. The elements of $\Delta$ are called  faces and the   maximal faces  under inclusion
are called facets. If $F_{1},\dots,F_{q}$ is a complete list of
the facets of $\Delta$, we  call $\{F_{1},\dots,F_{q}\}$,  denoted by $F(\Delta)$, the facet set of $\Delta$ and
usually write $\Delta$ as $\Delta=\langle F_{1},\dots,F_{q}\rangle$. The vertex set of $\Delta$ is denoted by $V(\Delta)$.
A subcollection of $\Delta$  is a simplicial complex whose facets are also facets of $\Delta$. For $U\subseteq V(\Delta)$, an
induced subcollection of $\Delta$ on $U$, denoted by $\Delta_{U}$, is the simplicial complex whose
facet set is $\{F\in F(\Delta)\mid F\subseteq U\}$ and vertex set is a subset of $U$.

If $F$ is a face of $\Delta=\langle F_{1},\dots,F_{q}\rangle$, we define the complement of F in $\Delta$ to be
$$F_{V}^{c}=V\setminus F\ \ \text{and}\ \ \Delta_{V}^{c}=\langle (F_{1})_{V}^{c},\dots,(F_{q})_{V}^{c}\rangle.$$
Note that if $U\subsetneq V(\Delta)$, we have that $\Delta_{U}^{c}=(\Delta_{U})_{U}^{c}$.
\end{Definition}

\vspace{3mm} To a squarefree monomial ideal $I$ in a polynomial ring $R=K[x_{1},\dots, x_{n}]$ over a field $K$, one can associate
two unique simplicial complexes  $\Delta(I)$ and $\mathcal{N}(I)$  on the vertex set labeled $\{x_{1},\dots, x_{n}\}$. Conversely given a simplicial complex $\Delta$ with vertices set $\{x_{1},\dots, x_{n}\}$, one can associate two unique squarefree
monomials $\mathcal{I}(\Delta)$  and  $\mathcal{N}(\Delta)$ in the polynomial ring $K[x_{1},\dots, x_{n}]$; these are all defined below
$$ \left. \begin{array}{ll}
\hspace{-2.5mm}\mbox{Facet complex of}\ I&\Delta(I)=\langle  \{x_{i_{1}},\dots, x_{i_{s}}\} | \prod\limits_{j=1}^{s}x_{i_{j}}\text{\  minimal generator of } I\rangle,\\
\hspace{-2.5mm}\mbox{Stanley-Reisner complex of}\ I&\mathcal{N}(I)=\langle \{x_{i_{1}},\dots, x_{i_{s}}\} |\prod\limits_{j=1}^{s}x_{i_{j}}\notin I \rangle,\\
\hspace{-2.5mm}\mbox{Facet ideal of}\ \Delta& \mathcal{I}(\Delta)=(\prod\limits_{j=1}^{s}x_{i_{j}}|\,\{x_{i_{1}},\dots, x_{i_{s}}\}\in
  F(\Delta)),\\
\hspace{-2.5mm}\mbox{Stanley-Reisner ideal of}\ \Delta& \mathcal{N}(\Delta)=(\prod\limits_{j=1}^{s}x_{i_{j}}|\,\{x_{i_{1}},\dots, x_{i_{s}}\}\notin \Delta).
\end{array}\right.
$$
Note that there is a one-to-one correspondence between monomial ideals and
simplicial complexes via each of these methods.

\vspace{3mm}For any homogeneous ideal $I$ in the polynomial ring $R=K[x_{1},\dots,x_{n}]$, there
exists a graded minimal free resolution of $R/I$
\vspace{3mm}
$$0\rightarrow \bigoplus\limits_{j}R(-j)^{\beta_{p,j}}\rightarrow \bigoplus\limits_{j}R(-j)^{\beta_{p-1,j}}\rightarrow \cdots\rightarrow \bigoplus\limits_{j}R(-j)^{\beta_{1,j}}\rightarrow R\rightarrow R/I\rightarrow 0,$$
where the maps are exact, $p\leq n$, and $R(-j)$ is the graded free module obtained by shifting
the degrees of $R$ by $j$. The number
$\beta_{i,j}$, the $(i,j)$-th graded Betti number of $R/I$,  is
an invariant of $R/I$ that equals the number of minimal generators of degree $j$ in the
$i$th syzygy module of $R/I$.
Of particular interest are the following invariants which measure the ※size§ of the minimal graded
free resolution of $R/I$.
The projective dimension of $R/I$, denoted $pd\,(R/I)$, is defined to be
$$pd\,(R/I):=\mbox{max}\{i\ |\ \beta_{i,j}\neq 0\}.$$
The regularity of $I$, denoted $reg\,(R/I)$, is defined by
$$reg\,(R/I):=\mbox{max}\{j-i\ |\ \beta_{i,j}\neq 0\}.$$

Let  $\Delta$ be a simplicial complex, we denote by $\mathcal{C}.(\Delta) $
 the reduced chain complex
and by $\widetilde{H}_{i}(\Delta)=Z_{i}(\Delta)/B_{i}(\Delta)$
 the $i$-th reduced homology groups of $\Delta$ with coefficients
in the field $K$. For computing the graded Betti numbers of the Stanley-Reisner ring of a
simplicial complex we use the following equivalent form of Hochster's formula.

\begin{theorem}\label{Thm1}(\cite[Theorem  2.8]{AF1})
Let $I$ be a pure squarefree monomial ideal in the polynomial ring $R=K[x_{1},\dots,x_{n}]$. Then the  graded Betti numbers of $R/I$ are given by
$$\beta_{i,j}(R/I)=\sum\limits_{\begin{array}{c}
\Gamma\subseteq \Delta(I),\\
|V(\Gamma)|=j\end{array}}dim_{K}\widetilde{H}_{i-2}(\Gamma_{V(\Gamma)}^{c}).$$
\end{theorem}

For the computation of the homology of certain subcomplexes as they appear
in the Hochster formula, we will repeatedly apply the following Mayer-Vietoris
sequence.

\vspace{3mm}
\begin{theorem}\label{Thm2}(\cite[Theorem  25.1]{M})
Let $\Delta$ be a simplicial complex, $\Delta_{1}$ and $\Delta_{2}$  subcomplexes of $\Delta$ such that $\Delta=\Delta_{1}\cup \Delta_{2}$. Then there is a long
exact sequence
$$\cdots\rightarrow \widetilde{H}_{i}(\Delta_{1}\cap \Delta_{2})\rightarrow \widetilde{H}_{i}(\Delta_{1})\oplus \widetilde{H}_{i}(\Delta_{2})\rightarrow \widetilde{H}_{i}(\Delta)\rightarrow \widetilde{H}_{i-1}(\Delta_{1}\cap \Delta_{2})\rightarrow \cdots$$
\end{theorem}

\vspace{3mm}\section{Path ideals  $I_{m,l}(C_{n})$ of a cycle graph  $C_{n}$}

In this section,  we will focus on   structures of path ideals of  a cycle graph $C_{n}$.

\begin{Definition} \label{def2}
Let $\Gamma$ be  a directed graph,  path ideal of $\Gamma$  of length $m$ is the  monomial ideal generated by some  paths of length $m$ in $\Gamma$, i.e.,
$$I_{m}(\Gamma)=(u_{1},\dots,u_{k}),\  \mbox{where}\  u_{1},\dots,u_{k}
\ \mbox{ are some paths  of length}\  m\ \mbox{in}\ \Gamma.$$
\end{Definition}
Note that if $u_{1},\dots,u_{k}$ are all paths of length $m$ in $\Gamma$, $I_{m}(\Gamma)=J_{m}(\Gamma)$, Futher, if $m=2$, then $I_{m}(\Gamma)$ is the edge ideal of graph $\Gamma$.

\vspace{3mm} In this paper, we consider the special case
that $\Gamma$  is an $n$-cycle  graph $C_{n}$, $u_{i}=\prod\limits_{j=1}^{m}x_{(i-1)l+j}$
 for any $1\leq i\leq k$ and  $k,l$ are positive integers such that $l<m$ and $k=\frac{n}{(l,n)}$, we denote this path ideal by $I_{m,l}(C_{n})$.

 \vspace{3mm}\begin{Example} \label{exm1}If $C_{4}$ is  a cycle with  vertex set $V=\{x_{1},x_{2},x_{3},x_{4}\}$, then we have
 $I_{3,2}(C_{4})=(x_{1}x_{2}x_{3},x_{3}x_{4}x_{1})$. If $C_{6}$ is a cycle with vertex set $V=\{x_{1},\dots,x_{6}\}$, then we have
 $I_{3,2}(C_{6})=(x_{1}x_{2}x_{3},x_{3}x_{4}x_{5},x_{5}x_{6}x_{1})$ and
 $I_{5,3}(C_{6})=(x_{1}x_{2}x_{3}x_{4}x_{5},x_{4}x_{5}x_{6}x_{1}x_{2})$.
\end{Example}

 \vspace{3mm}\begin{Note} \label{not1} Suppose that $l,l_{1},m,n$ are positive integers such that $l,l_{1}<m\leq n$. If $(l,n)=(l_{1},n)$,
  we have that $I_{m,l}(C_{n})=I_{m,l_{1}}(C_{n})$. In particular,
 we have that $I_{m,n-l}(C_{n})=I_{m,l}(C_{n})$. Hence we can suppose that $l$ is a divisor of $n$ such that $l\leq min\{m-1,\frac{n}{2}\}$.
 \end{Note}

 \vspace{3mm}\begin{Example} \label{exm2} Let $C_{12}$ be a cycle with vertex set $V=\{x_{1},\dots,x_{12}\}$, then we have
 $I_{11,1}(C_{12})=I_{11,5}(C_{12})=I_{11,7}(C_{12})=(\prod\limits_{j=1}^{11}x_{j},\dots,\prod\limits_{j=1}^{11}x_{(i-1)+j},\dots,\prod\limits_{j=1}^{11} x_{11+j})$,\\
$I_{11,2}(C_{12})=I_{11,10}(C_{12})=(\prod\limits_{j=1}^{11}x_{j},\prod\limits_{j=1}^{11}x_{2+j},\prod\limits_{j=1}^{11}x_{4+j},
\prod\limits_{j=1}^{11}x_{6+j},\prod\limits_{j=1}^{11}x_{8+j},\prod\limits_{j=1}^{11} x_{10+j})$,\\
$I_{11,3}(C_{12})=I_{11,9}(C_{12})=(\prod\limits_{j=1}^{11}x_{j},\prod\limits_{j=1}^{11}x_{3+j},\prod\limits_{j=1}^{11}x_{6+j},\prod\limits_{j=1}^{11} x_{9+j})$,
$I_{11,4}(C_{12})=I_{11,8}(C_{12})=(\prod\limits_{j=1}^{11}x_{j},\prod\limits_{j=1}^{11}x_{4+j},\prod\limits_{j=1}^{11} x_{8+j})$,
$I_{11,6}(C_{12})=(\prod\limits_{j=1}^{11}x_{j},\prod\limits_{j=1}^{11}x_{6+j})$,
 where we set $x_{a}=x_{i}$ when $a\equiv i$  mod $n$  $(1\leq i\leq n)$.
\end{Example}

\vspace{3mm}Throughout this paper, we will suppose that  $l,m,n$ are positive integers such that $2\leq m\leq n$ and $l\leq min\{m-1,\frac{n}{2}\}$ is a divisor of $n$. We also set $x_{j}=x_{i}$ when $j\equiv i$  mod $n$  $(1\leq i\leq n)$.

\vspace{3mm}\begin{Note}\label{note2} Let $C_{n}$ be an $n$-cycle graph, $l,m$ positive integers such that $2\leq m\leq n$, $l|n$ and $l\leq min\{m-1,\frac{n}{2}\}$ and let $k=\frac{n}{(l,n)}=\frac{n}{l}$.
We use the notation  $\Delta_{m,l}(C_{n})=\langle F_{1},\dots, F_{k}\rangle$ to denote  the   path complex, i.e., $F_{i}=\{x_{(i-1)l+1},\dots,x_{(i-1)l+m}\}$, $ i=1,\dots, k$, is the facet of $\Delta_{m,l}(C_{n})$,
 when $j\equiv i$  mod $n$  $(1\leq i\leq n)$, we set $x_{j}=x_{i}$.
 This labeling is called the
standard labeling of $\Delta_{m,l}(C_{n})$.

For any $1\leq i\leq k$, we have
$F_{i+1}\setminus F_{i}=\{x_{il+(m-l)+1},\dots,x_{il+m}\}$ and $F_{i}\setminus F_{i+1}=\{x_{(i-1)l+1},\dots,x_{(i-1)l+l}\}$,
it follows that $|F_{i}\setminus F_{i+1}|=|F_{i+1}\setminus F_{i}|=l$
 for all $1\leq i\leq k-1$.
\end{Note}

\vspace{3mm}
It is obvious that any induced subgraph of a cycle graph is a disjoint union of some
paths. Borrowing the terminology from  Jacques \cite{J}, we call the  path complex of
a line a ``run".

\begin{Definition} \label{def3} Given two positive  integers $l,m$ such that $l<m$, we define a run to be the path complex of a line graph. A run with $p$ facets is called a run of length $p$ and corresponds to $\Delta_{m,l}(L_{pl+(m-l)})$. Therefore a run of
length $p$ has $(p-1)l+m$ vertices.
\end{Definition}

\vspace{3mm}Set $supp(F_{i})=\{x_{j}: x_{j}\in F_{i}\}$ for $1\leq i\leq \frac{n}{(l,n)}$. By similar arguments as in \cite[Proposition 3.6]{AF1}, we can obtain that every proper induced subcollection of the path complex of a
cycle is a disjoint union of runs.

\begin{Remark} \label{rem1}Let $C_{n}$ be an $n$-cycle graph,  $l,m$ positive integers such that $2\leq m\leq n$, $l|n$ and $l\leq min\{m-1,\frac{n}{2}\}$. Let $\Gamma$ be
a proper induced connected subcollection of $\Delta_{m,l}(C_{n})$ on $U\subsetneq  X$. Then $\Gamma$ is  the line graph with the facet set $\{F_{i}:F_{i}\in \Delta_{m,l}(C_{n})\ \text{and} \ supp(F_{i})\subseteq U\}$.
\end{Remark}

\vspace{3mm}If $\Gamma$ and $\Lambda$ are two induced subcollections of $\Delta_{m,l}(C_{n})$
 composed of runs with equal lengths, one can easily see that $\Gamma$ and $\Lambda$  are homeomorphic as simplicial complexes
by a bijective map between their vertex sets. In particular, the  simplicial
complexes $\Gamma_{V(\Gamma)}^{c}$ and $\Lambda_{V(\Lambda)}^{c}$ are homeomorphic and have the same reduced homologies. Therefore, using Theorem \ref{Thm1} and
 Remark \ref{rem1}, all the information we need
to compute the Betti numbers of $\Delta_{m,l}(C_{n})$, or equivalently the homologies of induced
subcollections of $\Delta_{m,l}(C_{n})$, depend on the number and the lengths of the runs.

\vspace{3mm}\begin{Definition} \label{def4} Let $d$ be a positive integer,  a pure $d$-dimensional simplicial complex  $\Gamma=\langle F_{1},\dots, F_{s}\rangle$ be a disjoint union of runs of length $s_{1},\dots,s_{r}$. Then the sequence of positive
integers $s_{1},\dots,s_{r}$ is called a run sequence on $\mathcal{Y}=V(\Gamma)$, and we use the notation
$$E(s_{1},\dots,s_{r})=\Gamma_{\mathcal{Y}}^{c}=\langle (F_{1})_{\mathcal{Y}}^{c},\dots, (F_{s})_{\mathcal{Y}}^{c}\rangle.$$
\end{Definition}

\hspace{5mm}

\section{Reduced homologies for Betti numbers}

 Let $l,m,n$ be positive integers such that $2\leq m<n$, $l|n$ and $l\leq min\{m-1,\frac{n}{2}\}$, let $I=I_{m,l}(C_{n})$  be the path ideal of the the cycle $C_{n}$. By applying Theorem \ref{Thm1}, we see that to compute the Betti numbers  of
$R/I$, we need to compute the reduced homologies of complements of induced
subcollections of $\Delta_{m,l}(C_{n})$ which are disjoint unions of runs  by Remark \ref{rem1}. In this
section, we are  devote to complex homological calculations. The results  will allow
us to compute all Betti numbers of $R/I$  in the sections that follow. Firstly, we need the following lemma.

\vspace{3mm}
\begin{lemma}\label{lem1} (\cite[Lemma  4.1]{AF1}). Let $F_{1},\dots, F_{k}$ be subsets of a finite set $V$, where $k\geq 2$. Suppose that $\mathcal{E}=\langle (F_{1})_{V}^{c},\dots, (F_{k})_{V}^{c}\rangle$.
 Then for any $i$ we have

(1) Suppose $V\setminus \bigcup\limits_{j=2}^{k}F_{j}\neq\emptyset$. If $\mathcal{E}_{1}=\langle (F_{1})_{V}^{c}\rangle$ and $\mathcal{E}_{2}=\langle (F_{2})_{V}^{c},\dots, (F_{k})_{V}^{c}\rangle$,
 then
 $$ \widetilde{H}_{i}(\mathcal{E})=\widetilde{H}_{i-1}(\langle (F_{2})_{(V\setminus F_{1})}^{c},\dots,(F_{k})_{(V\setminus F_{1})}^{c}\rangle). $$

(2) If $F_{a}\subset F_{b}$ for some $a\neq b$, then $\mathcal{E}=\langle (F_{1})_{V}^{c},\dots,(\widehat{F_{b}})_{V}^{c},\dots,(F_{k})_{V}^{c}\rangle$.
\end{lemma}

The decomposition $\mathcal{E}=\mathcal{E}_{1}\cup \mathcal{E}_{2}$ which satisfies  conditions  above is called standard
decomposition of $\mathcal{E}$.

\vspace{3mm}Using the above lemma, we can obtain that

\begin{theorem} \label{Thm3}
Let $l,m,n$ be positive integers such that $2\leq m\leq n$, $l|n$ and $l\leq min\{m-1,\frac{n}{2}\}$,
let  $m\equiv s\ \mbox{mod}\ l$ with $0\leq s<l$, $t=\frac{m-s}{l}=1$ and $k=\frac{n}{l}$.  Let
$\Delta_{m,l}(C_{n})=\langle F_{1},\dots, F_{k}\rangle$ the path complex on the vertex set $V=\{x_{1},\dots,x_{n}\}$ with standard labeling. Suppose the connected
components of $\Delta_{m,l}(C_{n})$ are runs of lengths  $s_{1},\dots,s_{r}$, and $\mathcal{E}=E(s_{1},\dots,s_{r})$.  Then for all $i$, we have

(1) $$\widetilde{H}_{i}((\Delta_{m,l}(C_{n}))_{V}^{c})=\left\{\begin{array}{ll}
K,&\mbox{if}\ i=k-2\\
0,&\mbox{otherwise}.\\
\end{array}\right.$$

(2) $$\widetilde{H}_{i}(\mathcal{E})=\left\{\begin{array}{ll}
K,&\mbox{if}\ i=\sum\limits_{j=1}^{r}s_{j}-2\\
0,&\mbox{otherwise}.\\
\end{array}\right.$$
\end{theorem}
 \begin{pf} (1) can be shown by similar arguments as  (2), so we only prove  (2).

 According to Lemma \ref{lem1}, it is enough to show that this result is true for $r=1$.
There are two possible cases:

If $s_{1}=1$, then $\mathcal{E}=\langle (F_{1})_{V}^{c}\rangle=\{\emptyset\}$, and therefore
$$\widetilde{H}_{i}(\mathcal{E})=\left\{\begin{array}{ll}
K,&\mbox{if}\   i=-1\\
0,&\mbox{otherwise}.\\
\end{array}\right.$$

If $s_{1}\geq 2$, by repeated use
of Lemma \ref{lem1}, we get
$$\widetilde{H}_{i}(\mathcal{E})=\widetilde{H}_{i-(s_{1}-1)}(\{\emptyset\})=\left\{\begin{array}{ll}
K,&\mbox{if}\   i=s_{1}-2\\
0,&\mbox{otherwise}.\\
\end{array}\right.$$
The result follows.
\end{pf}

\vspace{3mm}\begin{Proposition} \label{prop2}
Let $l,m,n$ be positive integers such that $2\leq m\leq n$, $l|n$ and $l\leq min\{m-1,\frac{n}{2}\}$,
let  $m\equiv s\ \mbox{mod}\ l$ with $0\leq s<l$, $t=\frac{m-s}{l}\geq 2$ and $k=\frac{n}{l}$.  Let $\Delta_{m,l}(C_{n})=\langle F_{1},\dots, F_{k}\rangle$ be the path complex on the vertex set $V=\{x_{1},\dots,x_{n}\}$ with standard labeling. Suppose the connected
components of $\Delta_{m,l}(C_{n})$ are runs of lengths  $s_{1},\dots,s_{r}$, and $\mathcal{E}=E(s_{1},\dots,s_{r})$. Let $s_{j}=p_{j}(t+1)+d_{j}$
 where $p_{j}\geq 0$, $0\leq d_{j}\leq t$ and $1\leq j\leq r$. Then for all $i$, we have the
following conditions:

(1) If $s_{j}\geq t+2$, then $\widetilde{H}_{i}(\mathcal{E})\cong\widetilde{H}_{i-2p_{j}}(\langle (F_{1})_{W}^{c},\dots,(F_{s_{j-1}})_{W}^{c},(F_{p_{j}(t+1)+1})_{W}^{c},\dots,\\
(F_{p_{j}(t+1)+d_{j}})_{W}^{c},(F_{p_{j}(t+1)+d_{j}+1})_{W}^{c},\dots,(F_{s_{j+1}})_{W}^{c},\dots,(F_{k})_{W}^{c}\rangle),$
where $W=V\setminus A$ and $A$ is the  vertex set of  subcollection of the simplicial complex with
run of length  $s_{j}$, which composed of  the front $[(p_{j}-1)(t+1)+1]$  facets;

(2) If $s_{j}=2$ and $r\geq 2$, then $\widetilde{H}_{i}(\mathcal{E})=\widetilde{H}_{i-2}(\langle(F_{1})_{W}^{c},
\dots,(F_{s_{j-1}})_{W}^{c},(F_{s_{j+1}})_{W}^{c},\dots,\\
(F_{k})_{W}^{c}\rangle)$, where
$W$ is the vertex set of simplicial complex
$\langle  F_{1},\dots, F_{s_{j-1}},F_{s_{j+1}},\dots, F_{k}\rangle$;

(3) If $s_{j}=1$ and $r\geq 2$, then $\widetilde{H}_{i}(\mathcal{E})=\widetilde{H}_{i-1}(\langle
 (F_{1})_{W}^{c},\dots,(F_{s_{j-1}})_{W}^{c}, (F_{s_{j+1}})_{W}^{c},\dots,\\
 (F_{k})_{W}^{c}\rangle)$, where
$W$ is the vertex set of simplicial complex
$\langle F_{1},\dots, F_{s_{j-1}}, F_{s_{j+1}},\dots, F_{k}\rangle$;

(4) If  $d_{j}\neq 1,2$, then  $\widetilde{H}_{i}(\mathcal{E})=0$.
\end{Proposition}
 \begin{pf}
 Without loss of generality we can assume that $F_{1},\dots, F_{k}$  are ordered such that
$F_{1},\dots, F_{s_{j}}$
are the facets of the run of length $s_{j}$ and they have standard labeling
$F_{i}=\{x_{(i-1)l+1},\dots,x_{(i-1)l+m}\} $ for $1\leq i\leq s_{j}$.
We have $\mathcal{E}=\langle (F_{1})_{V}^{c},\dots, (F_{k})_{V}^{c}\rangle$.	Since $x_{1},\dots,x_{l}\in V\setminus \bigcup\limits_{j=2}^{k}F_{j}$, there is a standard decomposition
$$\mathcal{E}=\langle (F_{1})_{V}^{c}\rangle\cup \langle (F_{2})_{V}^{c},\dots, (F_{k})_{V}^{c}\rangle.$$
From Lemma \ref{lem1} (1), setting
$V'=V\setminus \{x_{1},\dots,x_{m}\}=\{x_{m+1},\dots,x_{n}\}$, we have that
$$\leqno{(4.1)\hspace{3.5cm}\widetilde{H}_{i}(\mathcal{E}) \cong\widetilde{H}_{i-1}(\langle (F_{2})_{V'}^{c},\dots,(F_{k})_{V'}^{c}\rangle).}$$
If  $s_{j}\geq t+2$, from  $(4.1)$, we have
$$\leqno{(4.2)\hspace{2.5cm} \widetilde{H}_{i}(\mathcal{E}) \cong\widetilde{H}_{i-1}(\langle\{x_{m+1},\dots,x_{m+l}\}_{V'}^{c},(F_{t+2})_{V'}^{c},\dots,(F_{k})_{V'}^{c}\rangle)}$$
and 	since $x_{m+1}\in V'\setminus \bigcup\limits_{j=t+2}^{k}F_{j}$, there is a standard decomposition
$$\langle\{x_{m+1},\dots,x_{m+l}\}_{V'}^{c}\rangle\cup \langle(F_{t+2})_{V'}^{c},\dots,(F_{s_{j}})_{V'}^{c}, \dots,(F_{k})_{V'}^{c}\rangle. $$
From (4.2) and Lemma \ref{lem1} $(1)$, with $V''=V\setminus \{x_{1},\dots,x_{m+l}\}=\{x_{m+l+1},\dots,x_{n}\}$, we have
\begin{eqnarray*}\widetilde{H}_{i}(\mathcal{E}) &\cong&\widetilde{H}_{i-2}(\langle(F_{t+2})_{V''}^{c},\dots,(F_{s_{j}})_{V''}^{c},\dots,(F_{k})_{V''}^{c}\rangle)\\
&=&\widetilde{H}_{i-2}(\langle\{x_{m+l+1},\dots,x_{m+(t+1)l}\}_{V''}^{c}, (F_{t+3})_{V''}^{c},\dots,(F_{s_{j}})_{V''}^{c},\dots,(F_{k})_{V''}^{c} \rangle).
\end{eqnarray*}
Again using the above arguments $(p_{j}-1)$ times to $s_{j}$, we can conclude that
$$\widetilde{H}_{i}(\mathcal{E})\cong
\widetilde{H}_{i-2p_{j}}(\langle (F_{p_{j}(t+1)+1})_{W}^{c},\dots,
(F_{p_{j}(t+1)+d_{j}})_{W}^{c},(F_{s_{j+1}})_{W}^{c},\dots,(F_{k})_{W}^{c}\rangle)$$
where $W=V\setminus \{x_{1},\dots,x_{(p_{j}-1)(t+1)l+m}\}$.
This settles Case $(1)$ of the proposition. Now suppose $1\leq s_{j}<t+2$. In this case by $(4.1)$ and
Lemma \ref{lem1}, we get that

$$\leqno{(4.3)}\hspace{1.3cm} \widetilde{H}_{i}(\mathcal{E}) \cong\widetilde{H}_{i-1}(\langle\{x_{m+1},\dots,x_{m+l}\}_{V'}^{c},(F_{s_{j}+1})_{V'}^{c},\dots,(F_{k})_{V'}^{c}\rangle) \ \ {for\ all}\ \ i.$$
\hspace{5mm}(a) If $s_{j}\geq 3$, since $x_{(s_{j}-1)l+(m-l)+1},\dots,x_{(s_{j}-1)l+m}\in V'\setminus (\bigcup\limits_{i=s_{j+1}}^{k}F_{j}\cup \{x_{m+1},\dots,x_{m+l}\})$, the simplicial complex
$\langle\{x_{m+1},\dots,x_{m+l}\}_{V'}^{c},(F_{s_{j}+1})_{V'}^{c},\dots,(F_{k})_{V'}^{c}\rangle$
is a cone, by Proposition $2.7$ of \cite{AF1} and $(4.3)$, we have $\widetilde{H}_{i}(\mathcal{E})=0$  for all $i$.

(b) If $s_{j}=2$ and $r\geq 2$, then $x_{m+1},\dots,x_{m+l}\in V'\setminus (\bigcup\limits_{i=s_{j+1}}^{k}F_{j})$,
we have  that $$\langle\{x_{m+1},\dots,x_{m+l}\}_{V'}^{c}\rangle\cup  \langle(F_{s_{j+1}})_{V'}^{c},\dots,(F_{k})_{V'}^{c}\rangle$$
is a standard decomposition. Thus  by Lemma \ref{lem1} and $(4.3)$, we have that
$$\widetilde{H}_{i}(\mathcal{E})\cong\widetilde{H}_{i-2}(\langle(F_{s_{j+1}})_{V''}^{c},\dots,(F_{k})_{V''}^{c}\rangle)\ \
 \text{for all}\  i.
$$
$V''$ is the vertex set of simplicial complex
$\langle F_{s_{j+1}},\dots, F_{k}\rangle$.
This settles Case $(2)$.

(c) If $s_{j}=1$ and $r\geq 2$, since $F_{1}\cap F_{h}=\emptyset$  for $1<h\leq k$, and from $(4.1)$ we obtain that  for all $i$,
 $\widetilde{H}_{i}(\mathcal{E})\cong\widetilde{H}_{i-1}(\langle (F_{s_{j+1}})_{V'}^{c},\dots,(F_{k})_{V'}^{c}\rangle)$
where $V'=V\setminus \{x_{1},\dots,x_{m}\}=\{x_{m+1},\dots,x_{n}\}$ is the vertex set of simplicial complex
$\langle F_{s_{j+1}},\dots, F_{k}\rangle$. This settles Case $(3)$.

 Now we prove Case $(4)$ by induction on $p_{j}$. If $p_{j}=0$, then $d_{j}=s_{j}\geq 1$.
 From above we know that $\widetilde{H}_{i}(\mathcal{E})=0$
 if  $s_{j}\geq 3$, and we are done. Now we can suppose $p_{j}\geq 1$. there are two possible cases:

(i) If $s_{j}<t+2$, then since $p_{j}\geq 1$, we must have $p_{j}=1$, $d_{j}=0$ and $s_{j}=t+1$. It was proved
above (under the case $s_{j}\geq 3$) that $\widetilde{H}_{i}(\mathcal{E})=0$.

(ii) If $s_{j}\geq t+2$, by Case $(1)$ we have
$$\widetilde{H}_{i}(\mathcal{E})\cong\widetilde{H}_{i-2p_{j}}(\langle (F_{p_{j}(t+1)+1})_{W}^{c},\dots,\\
(F_{p_{j}(t+1)+d_{j}})_{W}^{c},(F_{s_{j+1}})_{W}^{c},\dots,(F_{k})_{W}^{c}\rangle).
$$
From the above  equality and the case $p_{j}=0$, we know that $\widetilde{H}_{i}(\mathcal{E})=0$ for  $d_{j}\neq 1,2$. This proves $(4)$ and
we complete the proof.
\end{pf}

\vspace{3mm}
The following proposition  can be shown by similar arguments as in \cite[Proposition 4.3]{AF1}, we omit its proof

\vspace{3mm}\begin{Proposition} \label{prop3}
Let $l,m,n,k,t$ be  integers as in Proposition \ref{prop2}, let $\Delta_{m,l}(C_{n})=\langle F_{1},\dots, F_{k}\rangle$ be the path complex on the vertex set $V=\{x_{1},\dots,x_{n}\}$ with standard labeling.
Suppose that  $\alpha,\beta\geq 0$ are  integers and the connected
components of   $\Delta_{m,l}(C_{n})$ are runs of lengths  $p_{1}(t+1)+1,\dots,p_{\alpha}(t+1)+1$, $q_{1}(t+1)+2,\dots,q_{\beta}(t+1)+2$,
where $p_{i}, q_{j}\geq 0$ for $1\leq i\leq \alpha, 1\leq j\leq \beta$. Then
$$\widetilde{H}_{i}(\mathcal{E})=\left\{\begin{array}{ll}
K,&\mbox{if}\ i=2(P+Q)+2\beta+\alpha-2\\
0,&\mbox{otherwise},\\
\end{array}\right.$$
where $\mathcal{E}=E(p_{1}(t+1)+1,\dots,p_{\alpha}(t+1)+1,q_{1}(t+1)+2,\dots,q_{\beta}(t+1)+2)$, $P=\sum\limits_{i=1}^{\alpha}p_{i}$ and $Q=\sum\limits_{j=1}^{\beta}q_{j}$.
\end{Proposition}

\vspace{3mm}As a consequence of Propositions \ref{prop2} and \ref{prop3}, we obtain the following result:
\begin{corollary} \label{cor1}
Let $l,m,n,t$ be  integers as in Proposition \ref{prop2}.
Suppose that the  simplicial  complex $\Gamma=\langle F_{1},\dots, F_{p(t+1)+d}\rangle$ of  dimension $(m-1)$
 is a  run of length  $p(t+1)+d$ with standard labeling.  Let $\mathcal{E}=\langle (F_{1})_{V}^{c},\dots, (F_{p(t+1)+d})_{V}^{c}\rangle$.
Then $$\widetilde{H}_{i}(\mathcal{E})=\left\{\begin{array}{ll}
K,&\mbox{if}\  d=1, i=2p-1\\
K,&\mbox{if}\ d=2, i=2p\\
0,&\mbox{otherwise}.\\
\end{array}\right.$$
\end{corollary}
\begin{pf}  By Proposition \ref{prop2} $(4)$, if $d\neq 1,2$ the homology is zero. In the cases where $d=1,2$ the
result follows directly from Proposition  \ref{prop3}.
\end{pf}

\vspace{3mm}Before we prove the main results of this section, we need the following two lemmas.
\begin{lemma}\label{lem2} Let $l,m,n,k,s,t$ be  integers as in Proposition \ref{prop2}.  Suppose that
$\Delta_{m,l}(C_{n})$ $=\langle F_{1},\dots, F_{k}\rangle$ is the  path complex of the cycle  $C_{n}$ with standard labeling.
Let $b,c,e\in [n]$ be such that $c\leq t-1$ and $b+e+t-1<k$.
 Suppose $e=p(t+1)+d$, where $p\geq 0$, $0\leq d\leq t$. Set
$V=\{x_{(b-1)l+f+1},\dots, x_{(b+e+t-1)l}\}$  where $0\leq f<l$
 and $\mathcal{E}=\langle (F_{b}\setminus \{x_{(b-1)l+1},\dots,x_{(b-1)l+f}\})_{V}^{c},
 \dots,(F_{b+e-1})_{V}^{c},\{x_{(b+e+t-c-1)l+1},x_{(b+e+t-c-1)l+2},\\
 \dots,x_{(b+e+t-1)l}\}_{V}^{c}\rangle$.
Then for all $i$ we have
$$\widetilde{H}_{i}(\mathcal{E})=\left\{\begin{array}{ll}
K,&\mbox{if}\ d=1,i=2p\\
K,&\mbox{if}\ d=c+1,i=2p+1\\
0,&\mbox{otherwise}
\end{array}\right.$$
\end{lemma}
\begin{pf} Without loss of generality we can assume $b=1$ so that $V=\{x_{f+1},\dots,x_{(e+t)l}\}$
and
$\mathcal{E}=\langle (F_{1}\setminus \{x_{1},\dots,x_{f}\})_{V}^{c},\dots, (F_{e})_{V}^{c},\{x_{(e+t-c)l+1},x_{(e+t-c)l+2},\dots,x_{(e+t)l}\}_{V}^{c}\rangle$.\\
Since $x_{(e+t)l}\notin F_{h}$
for $1\leq h\leq e$, $\mathcal{E}$ has standard decomposition
$$\mathcal{E}=\langle (F_{1}\setminus \{x_{1},\dots,x_{f}\})_{V}^{c},\dots, (F_{e})_{V}^{c}\rangle\cup
\langle \{x_{(e+t-c)l+1},x_{(e+t-c)l+2},\dots,x_{(e+t)l}\}_{V}^{c}\rangle.$$
Set $V_{1}=V\setminus \{x_{(e+t-c)l+1},\dots,x_{(e+t)l}\}$,
by  Lemma \ref{lem1} and the fact that
$\{x_{(e-c)l+1},\\
\dots,x_{(e+t-c)l}\}\supset \{x_{(e-c+1)l+1},\dots,x_{(e+t-c)l}\}
\supset\cdots\supset \{x_{(e-1)l+1},\dots,x_{(e+t-c)l}\}$,
we obtain that
$$\leqno{(4.4)}\hspace{0.5cm}
\widetilde{H}_{i}(\mathcal{E})\cong \widetilde{H}_{i-1}((F_{1}\setminus \{x_{1},\dots,x_{f}\})_{V_{1}}^{c},\dots, (F_{e-c})_{V_{1}}^{c},
\{x_{(e-1)l+1},\dots, x_{(e+t-c)l}\}_{V_{1}}^{c}).
$$
 We prove our statement by induction on $|V|=(e+t)l-f=[p(t+1)+d+t]l-f$.  The base case is $|V|=(d+t)l-f$,
in which case $p=0$,  and $e=d\geq 1$.
There are two cases to be considered:

 (1) If $1\leq d\leq c$, then $e\leq c$, and so by $(4.4)$
 $$\widetilde{H}_{i}(\mathcal{E})\cong \widetilde{H}_{i-1}(\{x_{(e-1)l+1},\dots, x_{(e+t-c)l}\}_{V_{1}}^{c})$$
In this situation, the simplex $\{x_{(e-1)l+1},\dots, x_{(e+t-c)l}\}_{V_{1}}^{c}=\emptyset\Leftrightarrow (e-1)l+1\leq f+1 \Leftrightarrow d=e=1$, and hence we have
 $$\widetilde{H}_{i}(\mathcal{E})=\left\{\begin{array}{ll}
K,&\mbox{if}\ d=1,i=0,\\
0,&\mbox{otherwise}.
\end{array}\right.$$

 (2) If $d>c$, we use $(4.4)$ to note that since $x_{(e+t-c)l}\notin (F_{1}\setminus \{x_{1},\dots,x_{f}\})\cup\cdots\cup F_{e-c}$, we get that
 $\langle(F_{1}\setminus \{x_{1},\dots,x_{f}\})_{V_{1}}^{c},\dots, (F_{e-c})_{V_{1}}^{c}\rangle\cup
\langle\{x_{(e-1)l+1},\dots, x_{(e+t-c)l}\}_{V_{1}}^{c}\rangle$ is a standard decomposition.
 By Lemma \ref{lem1} and $(4.4)$ along with the fact that $e=d\leq t$, we find that if $V_{2}=V_{1}\setminus \{x_{(e-1)l+1},\dots,x_{(e+t-c)l}\}$, then
 \begin{eqnarray*}
 \widetilde{H}_{i}(\mathcal{E})&\cong & \widetilde{H}_{i-2}(\langle\{x_{f+1},\dots,x_{(e-1)l}\}_{V_{2}}^{c},
 (F_{2})_{V_{2}}^{c},\dots,
 \{x_{(e-c-1)l+1},\dots,x_{(e-1)l}\}_{V_{2}}^{c}\rangle)\\
 &\cong &\widetilde{H}_{i-2}(\langle\{x_{(e-c-1)l+1},\dots,x_{(e-1)l}\}_{V_{2}}^{c}\rangle).
 \end{eqnarray*}
 In this case, the simplex $\{x_{(e-c-1)l+1},\dots,x_{(e-1)l}\}_{V_{2}}^{c}=\emptyset\Leftrightarrow (e-c-1)l+1\leq f+1 \Leftrightarrow d=e=c+1$. Therefore
 $$\widetilde{H}_{i}(\mathcal{E})=\left\{\begin{array}{ll}
K,&\mbox{if}\ d=c+1,i=1\\
0,&\mbox{otherwise}.
\end{array}\right.$$
This settles the base case of the induction. Now suppose $|V|=(e+t)l-f>(d+t)l-f$ and the
theorem holds for all the cases where $|V|<(e+t)l-f$.
 Since $|V_{1}|=(e+t-c)l-f<|V|$,  we
use $(4.4)$ and  the induction hypothesis on $V_{1}$, now with the following
parameters: $c_{1}=t-c+1$, $e_{1}=e-c=p(t+1)+d-c$ and
$$d_{1}=\left\{\begin{array}{ll}
d-c,&\mbox{if}\ d\geq c,\\
d-c+t+1,&\mbox{if}\ d<c
\end{array}\right. \text{and}\ \  p_{1}=\left\{\begin{array}{ll}
p,&\mbox{if}\ d\geq c,\\
p-1,&\mbox{if}\ d<c
\end{array}\right.$$
Applying the induction hypothesis on $V_{1}$, we see that $\widetilde{H}_{i}(\mathcal{E})=0$, unless one of the
following cases happen, in which case  $\widetilde{H}_{i}(\mathcal{E})=K$.

1. $d_{1}=1$, and $i-1=2p_{1}$.

(a) When $d\geq c$, this means that $d=c+1$ and $i=2p+1$.

(b) When $d<c$, this means $d=c-t\geq 0$, which is not possible as we have assumed
$c\leq t-1$.

2. $d_{1}=c_{1}+1$,  and $i-1=2p_{1}+1$.

(a) When $d\geq c$, this means that $d-c=d_{1}=c_{1}+1$, and so $d=c+c_{1}+1=c+(t-c+1)+1=t+2$
which is not possible, as we have assumed $d\leq t$.

(b) When $d<c$, this means that $d-c+t+1=d_{1}=c_{1}+1=t-c+1+1$, and so $d=1$
and $i=2p_{1}+2=2p$.

We conclude that $\widetilde{H}_{i}(\mathcal{E})=K$ only when $d=1$ and $i=2p$, or $d=c+1$ and $i=2p+1$, and $\widetilde{H}_{i}(\mathcal{E})=0$ otherwise.
\end{pf}

\vspace{3mm}
\begin{lemma}\label{lem3} Let $l,m,n,k,t$ be  integers as in Proposition \ref{prop2}. Suppose that
$\Delta_{m,l}(C_{n})$ $=\langle F_{1},\dots, F_{k}\rangle$ is the  path complex of the cycle  $C_{n}$ on vertex set $V=\{x_{1},\dots,x_{n}\}$ with standard labeling. We can write $k$ as  $k=p(t+1)+d$, where $p\geq 1$, $0\leq d\leq t$.
Consider the following simplicial complexes:
\begin{eqnarray*}(4.5)\hspace{1.3cm}&E_{0}&=\langle (F_{1})_{V}^{c},(F_{2})_{V}^{c}\dots, (F_{k-t+1})_{V}^{c}\rangle=E(k-t+1),\\
&E_{a}&=E_{a-1}\cup\langle (F_{k-a+1})_{V}^{c}\rangle\ \ \text{for}\ \ 1\leq a\leq t-1.
 \end{eqnarray*}
Then for all $0\leq a\leq t-2$, we have
$$\widetilde{H}_{i}(E_{a}\cap \langle(F_{k-a})_{V}^{c}\rangle)=\left\{\begin{array}{ll}
K,&\mbox{if}\  d=0, i=2p-3\\
K,&\mbox{if}\   d=t-a-1,i=2p-2\\
0,&\mbox{otherwise}.\\
\end{array}\right.
$$
\end{lemma}

\begin{pf} Setting $V'=V\setminus F_{k-a}=\{x_{(k-a-1)l+m+1-n},x_{(k-a-1)l+m+2-n},\dots$, $x_{(k-a-1)l}\}$, where if $j\equiv i$  mod $n$ $(1\leq i\leq n)$, we set $x_{j}=x_{i}$.
We can write $$E_{a}\cap \langle(F_{k-a})_{V}^{c}\rangle=\langle(F_{1})_{V'}^{c},\dots, (F_{k-t+1})_{V'}^{c},
(F_{k-a+1})_{V'}^{c},\dots, (F_{k})_{V'}^{c}\rangle.$$
We now compute the 	$(F_{h})_{V'}^{c}=\{x_{(h-1)l+1},\dots, x_{(h-1)l+m}\}_{V'}^{c}$
appearing above.

(1)  If $h$ is chosen as $k-a+1,\dots,k,1,\dots, t-a$, then we have $(F_{k-a+1})_{V'}^{c}\supset\cdots \supset (F_{k})_{V'}^{c}\supset (F_{1})_{V'}^{c}\supset\cdots \supset (F_{t-a})_{V'}^{c}$ and $(F_{k-a+1})_{V'}^{c}=\{x_{(k-a-1)l+m+1-n},\dots,x_{(k-a)l+m-n}\}_{V'}^{c}$.

\vspace{3mm}(2)If $h$ is chosen as $k-t-a,\dots,k-t+1$, then we have
$(F_{k-t-a})_{V'}^{c}\subset \dots\subset (F_{k-t+1})_{V'}^{c}$
  and $(F_{k-t+1})_{V'}^{c}=\{x_{(k-t)l+1},\dots,x_{(k-a-1)l}\}_{V'}^{c}$.

\vspace{3mm}(3) If $t-a+1\leq h\leq k-t-a-1$,
then $$(F_{h})_{V'}^{c}=\{x_{(k-a-1)l+m+1-n},\dots,x_{(h-1)l}, x_{(h-1)l+m+1},\dots,x_{(k-a-1)l}\}. $$

From the above observations,  we see that

 \begin{eqnarray*}(4.6)\hspace{0.5cm}E_{a}\cap \langle(F_{k-a})_{V}^{c}\rangle\!&\!=\!&\!
\langle\{x_{(k-a-1)l+m+1-n},\dots,x_{(k-a)l+m-n}\},F_{t-a+1},\dots,F_{k-t-a-1},\\
& &\{x_{(k-t)l+1},\dots,x_{(k-a-1)l}\}\rangle_{V'}^{c}.
\end{eqnarray*}

We now consider the following scenarios.

 Suppose $p=1$. In this situation, $k=(t+1)+d\leq 2t+1$ which implies that $k-t-a-1\leq t-a$.  Hence $(4.6)$ becomes
\begin{eqnarray*}(4.7)\hspace{1.5cm}E_{a}\!\cap\! \langle(F_{k-a})_{V}^{c}\rangle\!\!&\!\!=\!\!&\!\!
\langle\{x_{(k-a-1)l+m+1-n},\dots,x_{(k-a)l+m-n}\},
\{x_{(k-t)l+1},\dots,\\
& &x_{(k-a-1)l}\}\rangle_{V'}^{c}.
\end{eqnarray*}
Set $A=(k-a-1)l+m+1-n$, $B=(k-a)l+m-n$, $C=(k-t)l+1$ and $D=(k-a-1)l$. By the assumption that $m=tl+s$ and $n=kl=[(t+1)+d]l$,  we can obtain that $B-D=(k-a)l+m-n-(k-a-1)l=s-dl$ and $A-C=(k-a-1)l+m+1-n-[(k-t)l+1]=(t-a-2-d)l+s$.

(a) If $d\leq t-a-2$, then $A\geq C$. We  consider the following two cases:

 (i) If $s\leq dl$, then $B\leq D$. Hence
 $\{x_{(k-a-1)l+m+1-n},\dots,x_{(k-a)l+m-n}\}\subseteq\{x_{(k-t-1)l+1},\\
 \dots,x_{(k-a-1)l}\}$, which means that $(4.7)$ becomes
 $$E_{a}\cap \langle(F_{k-a})_{V}^{c}\rangle=\langle\{x_{(k-a-1)l+m+1-n},\dots,x_{(k-a)l+m-n}\} \rangle_{V'}^{c}.$$
Also note that  the simplex $\{x_{(k-a-1)l+m+1-n},\dots,x_{(k-a)l+m-n}\}_{V_{'}}^{c}=\emptyset\Leftrightarrow (k-a)l+m-n\geq(k-a-1)l\Leftrightarrow (t+1)l+dl\leq tl+l+s
\Leftrightarrow d=s=0$.
 It follows that $$\widetilde{H}_{i}(E_{a}\cap \langle(F_{k-a})_{V}^{c}\rangle)\cong\left\{\begin{array}{ll}
K,&\mbox{if}\  d=s=0, i=-1\\
0,&\mbox{otherwise}.
\end{array}\right.$$

(ii) If  $s>dl$, from the fact that $s<l$, we can obtain that $d=0$, $s>0$, $l\geq 2$ and $B>D$. Therefore
$A-C=(t-a-2-d)l+s=(t-a-2)l+s>0$, $A-D=(k-a-1)l+m+1-n-[(k-a-1)l]=m+1-n=s+1-l\leq 0$.
We can apply Lemma \ref{lem1},with $V''=V'\setminus \{x_{(k-t)l+1},\dots,x_{(k-a-1)l}\}=\{x_{(k-a-1)l+m+1-n},\dots,x_{(k-t)l}\}$ to find that for all $i$
  $$\widetilde{H}_{i}(E_{a}\cap \langle(F_{k-a})_{V}^{c}\rangle)=\widetilde{H}_{i-1}(\langle\{x_{(k-a-1)l+m+1-n},\dots,x_{(k-a)l+m-n}\} \rangle_{V''}^{c}).
 $$
In this case, $\langle\{x_{(k-a-1)l+m+1-n},\dots,x_{(k-a)l+m-n}\} \rangle_{V''}^{c}=\emptyset $,
it implies that $$\widetilde{H}_{i}(E_{a}\cap \langle(F_{k-a})_{V}^{c}\rangle)\cong\left\{\begin{array}{ll}
K,&\mbox{if}\  i=-1\\
0,&\mbox{otherwise}.
\end{array}\right.$$
Therefore, when $d\leq t-a-2$, we have
$$\widetilde{H}_{i}(E_{a}\cap \langle(F_{k-a})_{V}^{c}\rangle)\cong\left\{\begin{array}{ll}
K,&\mbox{if}\ d=0, i=-1\\
0,&\mbox{otherwise}.
\end{array}\right.$$

(b) If $d=t-a-1$, then $A-C=(t-a-2-d)l+s<0$ and $B-D=s-dl<0$.
Again using Lemma \ref{lem1}, with $V''=V'\setminus \{x_{(k-t)l+1},\dots,x_{(k-a-1)l}\}=\{x_{(k-a-1)l+m+1-n},\dots,x_{(k-t)l}\}$, one has
  $$\widetilde{H}_{i}(E_{a}\cap \langle(F_{k-a})_{V}^{c}\rangle)=\widetilde{H}_{i-1}(\langle\{x_{(k-a-1)l+m+1-n},\dots,x_{(k-a)l+m-n}\} \rangle_{V''}^{c}).
 $$
  In this case  $\langle\{x_{(k-a-1)l+m+1-n},\dots,x_{(k-a)l+m-n}\} \rangle_{V''}^{c}=\emptyset $, and
 it follows that $$\widetilde{H}_{i}(E_{a}\cap \langle(F_{k-a})_{V}^{c}\rangle)\cong\left\{\begin{array}{ll}
K,&\mbox{if}\  d=t-a-1,i=-1\\
0,&\mbox{otherwise}.
\end{array}\right.$$

(c) If $d>t-a-1$, then $A-C<0$, $B-C=(k-a)l+m-n-[(k-t)l+1]=(t-a)l+m-n-1=(t-a-1-d)l+s-1<0$.
We can apply Lemma \ref{lem1}, setting $V''=V'\setminus \{x_{(k-t)l+1},\dots,x_{(k-a-1)l}\}=\{x_{(k-a-1)l+m+1-n},\dots,x_{(k-t)l}\}$, to  obtain that
   $$\widetilde{H}_{i}(E_{a}\cap \langle(F_{k-a})_{V}^{c}\rangle)=\widetilde{H}_{i-1}(\langle\{x_{(k-a-1)l+m+1-n},\dots,x_{(k-a)l+m-n}\} \rangle_{V''}^{c}).
 $$
In this case $\langle\{x_{(k-a-1)l+m+1-n},\dots,x_{(k-a)l+m-n}\} \rangle_{V''}^{c}$ is a simplex.
Therefore, we can conclude that $\widetilde{H}_{i}(E_{a}\cap \langle(F_{k-a})_{V}^{c}\rangle)=0$.

Now suppose that $p\geq 2$. In this case it is easy to see that $A<(t-a)l+1$, $(k-t-a-2)l+1$,
$(k-t)l+1$. Therefore, we can apply Lemma \ref{lem1} (1) with $V''=V'\setminus \{x_{(k-a-1)l+m+1-n},\dots,x_{(k-a)l+m-n}\}=\{x_{(k-a)l+m+1-n},\dots,x_{(k-a-1)l}\}$ to
$(4.6)$ to conclude that for all $i$
\begin{eqnarray*}
\widetilde{H}_{i}(E_{a}\cap \langle(F_{k-a})_{V}^{c}\rangle)&=&\widetilde{H}_{i-1}(\langle (F_{t-a+1}\setminus \{x_{(t-a)l+1},\dots,x_{(k-a)l+m-n}\}),F_{t-a+2},\dots,\\
& &F_{k-t-a-1},\{x_{(k-t)l+1},\dots,x_{(k-a-1)l}\}\rangle_{V''}^{c}).
 \end{eqnarray*}
Note that $|V''|=[(p-1)(t+1)+d]l-s$. We can use  Lemma \ref{lem2} with values $f=s$, $b=t-a+1$, $c=t-a$, and $e=k-2t-1=(p-2)(t+1)+d+1$ to conclude that
$$\widetilde{H}_{i}(E_{a}\cap \langle(F_{k-a})_{V}^{c}\rangle)=\left\{\begin{array}{ll}
K,&\mbox{if}\  d=0, i=2p-3;\\
K,&\mbox{if}\  d=t-a-1, i=2p-2;\\
0,&\mbox{otherwise}.\\
\end{array}\right.
$$
\end{pf}

\vspace{3mm}We end this section with the calculation of the homology of the complement
of the  path complex of a cycle; this will give us the top degree Betti numbers
of the path ideal of a cycle. The following are the main results of this section.

\begin{theorem}\label{Thm4}Let $l,m,n,k,t$ be  integers as in Proposition \ref{prop2}.  Suppose that
$\Delta_{m,l}(C_{n})$ $=\langle F_{1},\dots, F_{k}\rangle$ is the  path complex of the cycle  $C_{n}$ on vertex set $V$ with standard labeling.
 We can write $k$ as  $k=p(t+1)+d$, where $p\geq 0$, $0\leq d\leq t$. Then for all $i$
$$\widetilde{H}_{i}((\Delta_{m,l}(C_{n}))_{V}^{c})=\left\{\begin{array}{ll}
K^{t},&\mbox{if}\  d=0,i=2p-2,p>0;\\
K,&\mbox{if}\ d\neq 0, i=2p-1;\\
0,&\mbox{otherwise}.\\
\end{array}\right.$$
\end{theorem}
\begin{pf}  By the previous assumptions that $m\leq n$, $t=\frac{m-s}{l}$ and $k=\frac{n}{l}$, we have that $t\leq k$.
If $p=0$, then $k=d\leq t$. This shows that $k=d=t$. In this situation $m=n$ and our claim is obvious, so we assume that $p\geq 1$, and therefore, $k\geq t+1$.

Considering the simplicial complexes in $(4.5)$,  we obtain that $(\Delta_{m,l}(C_{n}))_{V}^{c}=E_{t-1}$.
First, one has  $$k-t+1=p(t+1)+d-t+1=\left\{\begin{array}{ll}
p(t+1)+1,&\mbox{if}\  d=t\\
p(t+1),&\mbox{if}\  d=t-1\\
(p-1)(t+1)+d+2,&\mbox{if}\ d\leq t-2.\\
\end{array}\right.$$
Using the fact that
$E_{0}=E(k-t+1)$ and Corollary \ref{cor1}, we obtain that
$$ (4.8)\hspace{1.5cm}\widetilde{H}_{i}(E_{0})=\left\{\begin{array}{ll}
K,&\mbox{if}\  d=0, i=2p-2\\
K,&\mbox{if}\ d=t, i=2p-1\\
0,&\mbox{otherwise}.\\
\end{array}\right.$$
To calculate the homologies of $E_{t-1}$, we shall repeatedly apply the Mayer-Vietoris sequence as follows.
For a fixed $1\leq a\leq t-1$, since  $\langle(F_{k-a+1})_{V}^{c}\rangle$ is a simplex,  we have $\widetilde{H}_{i}(\langle(F_{k-a+1})_{V}^{c}\rangle)=0$ for all $i$.  Hence there exists the following exact
sequence:

$$\leqno{(4.9)\hspace{0.5cm}\widetilde{H}_{i}(E_{a-1}\cap \langle(F_{k-a+1})_{V}^{c}\rangle)\to \widetilde{H}_{i}(E_{a-1})\to \widetilde{H}_{i}(E_{a})\to \widetilde{H}_{i-1}(E_{a-1}\cap \langle(F_{k-a+1})_{V}^{c}\rangle)}.$$
We distinguish the following  cases:

(1) If $0<d<t$, then by Lemma \ref{lem3}, we know that
$$\widetilde{H}_{i}(E_{a-1}\cap \langle(F_{k-a+1})_{V}^{c}\rangle)=\left\{\begin{array}{ll}
K,&\mbox{if}\  a=t-d, i=2p-2\\
0,&\mbox{otherwise}.\\
\end{array}\right.$$
 We apply this observation and $(4.8)$
to the exact sequence  $(4.9)$ to see that for all $i$
$$\widetilde{H}_{i}(E_{a})=\widetilde{H}_{i}(E_{a-1})=\cdots=\widetilde{H}_{i}(E_{0})=0\ \ \mbox{for\ any}\ \ 1\leq a\leq t-d-1.$$
Once again we use  $(4.9)$  to obtain that for all $i$
\begin{eqnarray*}\widetilde{H}_{i}(E_{a})&\cong&\left\{\begin{array}{ll}
0,&\mbox{if}\  1\leq a\leq t-d-1\\
\widetilde{H}_{i-1}(E_{a-1}\cap \langle(F_{k-a+1})_{V}^{c}\rangle),&\mbox{if}\  a=t-d\\
\widetilde{H}_{i}(E_{t-d}),&\mbox{if}\  t-d<a\leq t-1\\
\end{array}\right.\\
&=&\left\{\begin{array}{ll}
K,&\mbox{if}\  a\geq t-d, i=2p-1\\
0,&\mbox{otherwise}.\\
\end{array}\right.
\end{eqnarray*}
Therefore, we can conclude that in this case
$$\widetilde{H}_{i}(E_{t-1})=\left\{\begin{array}{ll}
K,&\mbox{if}\ i=2p-1\\
0,&\mbox{otherwise}.\\
\end{array}\right.$$

(2) If $d=t$, then by Lemma \ref{lem3} we know that $\widetilde{H}_{i}(E_{a-1}\cap \langle(F_{k-a+1})_{V}^{c}\rangle)=0$ for all $i$ and all $1\leq a\leq t-1$. Applying  this fact along with $(4.8)$ to the sequence in $(4.9)$, we see that for $1\leq a\leq t-1$,
$$\widetilde{H}_{i}(E_{a})\cong\widetilde{H}_{i}(E_{0})=\left\{\begin{array}{ll}
K,&\mbox{if}\  i=2p-1\\
0,&\mbox{otherwise}.\\
\end{array}\right.$$

(3) If $d=0$,  then by Lemma \ref{lem3}, we know that
$$\widetilde{H}_{i}(E_{a-1}\cap \langle(F_{k-a+1})_{V}^{c}\rangle)=\left\{\begin{array}{ll}
K,&\mbox{if}\ i=2p-3\\
0,&\mbox{otherwise}.\\
\end{array}\right.$$
On the other hand, by $(4.8)$, we know $$\widetilde{H}_{i}(E_{0})=\left\{\begin{array}{ll}
K,&\mbox{if}\  i=2p-2\\
0,&\mbox{otherwise}.\\
\end{array}\right.$$
Applying
these facts to $(4.9)$, we see that
$$\widetilde{H}_{i}(E_{a})=\widetilde{H}_{i}(E_{a-1})=\cdots=\widetilde{H}_{i}(E_{0})=0\ \ \mbox{for}\ \ i\neq 2p-2.$$
When $i=2p-2$, the sequence $(4.9)$ produces an exact sequence

$$0\to \widetilde{H}_{2p-2}(E_{0})\to \widetilde{H}_{2p-2}(E_{1})\to \widetilde{H}_{2p-3}(E_{0}\cap \langle(F_{k})_{V}^{c}\rangle)\to 0.$$
Since $\widetilde{H}_{2p-2}(E_{0})=\widetilde{H}_{2p-3}(E_{0}\cap \langle(F_{k})_{V}^{c}\rangle)=K$, we have
$$\widetilde{H}_{i}(E_{1})=\left\{\begin{array}{ll}
K^{2},&\mbox{if}\  i=2p-2\\
0,&\mbox{otherwise}.\\
\end{array}\right.$$
We repeat the above method, recursively, for values $a=2,\dots,t-1$
$$0\to \widetilde{H}_{2p-2}(E_{a-1})\to \widetilde{H}_{2p-2}(E_{a})\to \widetilde{H}_{2p-3}(E_{a-1}\cap \langle(F_{k-a+1})_{V}^{c}\rangle)\to 0$$
and conclude that for $1\leq a\leq t-1$
$$\widetilde{H}_{i}(E_{a})=\left\{\begin{array}{ll}
K^{a+1},&\mbox{if}\  i=2p-2\\
0,&\mbox{otherwise}.\\
\end{array}\right.$$
We put this all together
$$\widetilde{H}_{i}(E_{t-1})=\left\{\begin{array}{ll}
K^{t},&\mbox{if}\  d=0, i=2p-2,p>0\\
K,&\mbox{if}\ d\neq 0,i=2p-2\\
0,&\mbox{otherwise}\\
\end{array}\right.$$
and this completes the proof.
\end{pf}

\vspace{5mm}

\vspace{3mm}\section{Graded Betti numbers of path ideals}

In this section,  we compute the top degree Betti numbers of path ideals  by applying the homological calculations from the previous
section.

\vspace{3mm}By Theorem \ref{Thm1} and Theorem \ref{Thm2}, we can compute the Betti numbers of $I_{m,l}(C_{n})$  of
degree $n$.
\begin{theorem}\label{Thm5}Let $l,m,n,k,t $ be  integers as in Proposition \ref{prop2}.  Suppose that
$C_{n}$ is an $n$-cycle graph. We write $k$ as  $k=p(t+1)+d$, where $p\geq 0$, $0\leq d\leq t$. Then
$$\beta_{i,n}(R/I_{m,l}(C_{n}))=\left\{\begin{array}{ll}
t,&\mbox{if}\  d=0, i=2p;\\
1,&\mbox{if}\ d\neq 0, i=2p+1;\\
0,&\mbox{otherwise}.\\
\end{array}\right.$$
\end{theorem}
\begin{pf}
 By Theorem \ref{Thm1},
$\beta_{i,n}(R/I_{m,l}(C_{n}))=dim_{K}\widetilde{H}_{i-2}((\Delta_{m,l}(C_{n}))_{V}^{c})$, and
the result now follows directly from Theorem \ref{Thm4}.
\end{pf}

\vspace{3mm}From Hochster's formula we see that computing Betti numbers of degree less
than n comes down to counting induced subcollections of certain kinds.

\begin{theorem}\label{Thm7} Let $l,m,n,k,t $ be  integers as in Proposition \ref{prop2}.  Suppose that
$\Delta_{m,l}(C_{n})$ $=\langle F_{1},\dots, F_{k}\rangle$ is the  path complex of the cycle  $C_{n}$ with standard labeling.
Let $\Lambda$ be an induced subcollection  of
$\Delta_{m,l}(C_{n})$ and $I=\mathcal{I}(\Lambda)$   the facet ideal of $\Lambda$. Let $F_{1},\dots, F_{\alpha}$,
$F_{\alpha+1},\dots, F_{\alpha+\beta}$ be all connected
components of $\Lambda$ with disjoint runs of lengths $p_{1}(t+1)+1,\dots,p_{\alpha}(t+1)+1$, $q_{1}(t+1)+2,\dots,q_{\beta}(t+1)+2$ where
 $\alpha,\beta, p_{u}, q_{v}\geq 0$ are  integers.
Assume that $i$ and $j$ are integers with $i\leq j<n$. Then the graded Betti number
$\beta_{i,j}(R/I)$
 is the number of subcollections of $\Lambda$ satisfying the
conditions $j=[(P+Q)(t+1)+\beta]l+m(\alpha+\beta)$, $i=2(P+Q)+2\beta+\alpha$ where $P=\sum\limits_{i=1}^{\alpha}p_{i}$ and $Q=\sum\limits_{j=1}^{\beta}q_{j}$.
\end{theorem}
\begin{pf} Since $\Delta(I)=\Lambda$, from Theorem \ref{Thm1}, we have
$$\beta_{i,j}(R/I)=\sum\limits_{\begin{array}{c}
\Gamma\subseteq \Lambda,\\
|V|=j\end{array}}dim_{k}\widetilde{H}_{i-2}(\Gamma_{V}^{c}).$$
where $V$ is the vertex set of $\Gamma$ and the sum is taken over induced subcollections $\Gamma$ of
$\Delta_{m,l}(C_{n})$.

Each induced subcollection of $\Lambda$ is clearly an induced subcollection of $\Delta_{m,l}(C_{n})$
and can therefore be written as a disjoint union of runs. So from Proposition \ref{prop2} $(4)$,
we can conclude the only $\Gamma$ whose complements have nonzero homology are those composed of disjoint runs of lengths $p_{1}(t+1)+1,\dots,p_{\alpha}(t+1)+1$, $q_{1}(t+1)+2,\dots,q_{\beta}(t+1)+2$. Such subcollections have $j$
vertices, where by Definition \ref{def3}
\begin{eqnarray*}(5.1)\hspace{0.5cm}j&=&p_{1}(t+1)l+m+\cdots+p_{\alpha}(t+1)l+m\\
&+&(q_{1}(t+1)+1)+m+\cdots+(q_{\beta}(t+1)+1)+m\\
&=&[(P+Q)(t+1)+\beta]l+m(\alpha+\beta).
\end{eqnarray*}
Hence
$$\Gamma_{V}^{c}=E(p_{1}(t+1)+1,\dots,p_{\alpha}(t+1)+1,q_{1}(t+1)+2,\dots,q_{\beta}(t+1)+2)$$
and by Proposition \ref{prop3}, we have
$$ (5.2)\hspace{1.5cm}dim_K(\widetilde{H}_{i-2}(\Gamma_{V}^{c}))=\left\{\begin{array}{ll}
1,&\mbox{if}\ i=2(P+Q)+2\beta+\alpha\\
0,&\mbox{otherwise}\\
\end{array}\right.$$
From $(5.1)$ and $(5.2)$,  we see that each induced subcollection $\Gamma$ with runs of lengths $p_{u}(t+1)+1$ or
 $q_{v}(t+1)+2$  contributes 1 unit to
$\beta_{i,j}(R/I)$ if and only if
$$
j=[(P+Q)(t+1)+\beta]l+m(\alpha+\beta),\
i=2(P+Q)+2\beta+\alpha.$$
\end{pf}

\begin{theorem}\label{Thm8}Let $i,j,l,m,n,k,t,\alpha,\beta,P,Q$ be  integers as in Theorem \ref{Thm4}.
  Suppose that $\Delta_{m,l}(C_{n})$ $=\langle F_{1},\dots, F_{k}\rangle$ is the  path complex of the cycle  $C_{n}$ with standard labeling.
 We can write $k$ as $k=p(t+1)+d$ with $p\geq 0$, $0\leq d\leq t$. Assume that $i$ and $j$ are integers with $i\leq j<n$. If
$\beta_{i,j}(R/I_{m,l}(C_{n}))\neq 0$, we have

(1) $j\leq mi$;

(2) If $d=0$, then $i<2p$;

(3) If $d\neq 0$, then $i\leq 2p+1$;

(4) $j-i\leq \left\{\begin{array}{ll}n-2p& \mbox{if}\ \ d=0\\
n-2p-2& \mbox{if}\ \ d\neq 0.
\end{array}\right.$
\end{theorem}
\begin{pf} By \cite[Theorem 3.3.4]{J}, we know that
$\beta_{i,j}(R/I_{m,l}(C_{n}))=0\ \ \mbox{for all}\ j>mi.$
This settles Case (1).
By  Theorem \ref{Thm4}, we know
$\beta_{i,j}(R/I_{m,l}(C_{n}))$ is equal to the number of all subcollections of $\Delta_{m,l}(C_{n})$, which are
 composed of disjoint runs of lengths $p_{1}(t+1)+1,\dots,p_{\alpha}(t+1)+1$, $q_{1}(t+1)+2,\dots,q_{\beta}(t+1)+2$ and
these length satisfying the
conditions $j=[(P+Q)(t+1)+\beta]l+m(\alpha+\beta)$, $i=2(P+Q)+2\beta+\alpha$.
Therefore, using assumption $m=tl+s$, we can obtain
\begin{eqnarray*}(5.3)\hspace{0.5cm}j-li&=&[(P+Q)(t+1)+\beta]l+m(\alpha+\beta)-l[2(P+Q)+2\beta+\alpha]\\
&=&(P+Q)(t-1)l+(m-l)(\alpha+\beta)\\
&=&[(P+Q)+(\alpha+\beta)](t-1)l+s(\alpha+\beta),
\end{eqnarray*}
\begin{eqnarray*}\hspace{1.5cm}tli-j&=&tl[2(P+Q)+2\beta+\alpha]-[(P+Q)(t+1)+\beta]l-m(\alpha+\beta)\\
&=&(P+Q+\beta)(t-1)l+(tl-m)(\alpha+\beta)\\
&=&(P+Q+\beta)(t-1)l-s(\alpha+\beta).
\end{eqnarray*}
It follows that $j-li\geq tli-j$.
It follows that $$(t+1)il\leq 2j\ \Rightarrow\ i\leq \frac{2j}{(t+1)l}<\frac{2n}{(t+1)l}.$$
As $n=kl=[p(t+1)+d]l$, we can obtain that $i<\frac{2[p(t+1)+d]}{(t+1)}$.
Therefore, if $d=0$ it follows that $i<2p$, and if $d\neq 0$ it follows that $i\leq 2p+1$.
This settles Cases (2) and (3).

On the one hand, since $\Delta_{m,l}(C_{n})$
 has $k$ facets and since there must be at least
$t$ facets between every two runs in $\Delta_{m,l}(C_{n})$, we have
$$k\geq (P+Q)(t+1)+\alpha+2\beta+t(\alpha+\beta),$$
it follows that $P+Q+\alpha+\beta\leq \frac{k-\beta}{t+1}$.

On the other hand, by assumption $m=tl+s$, $n=kl$ and  $k=p(t+1)+d$, we have
\begin{eqnarray*}\hspace{1.5cm}j-i&=&[(P+Q)(t+1)+\beta]l+m(\alpha+\beta)-[2(P+Q)+2\beta+\alpha]\\
&=&(P+Q+\alpha+\beta)[(t+1)l-2]+\beta l+\alpha+(s-l)(\alpha+\beta)\\
&\leq &\frac{k-\beta}{t+1}(t+1)l-\frac{2(k-\beta)}{t+1}+\beta l+\alpha+(s-l)(\alpha+\beta)\\
&=&kl-\frac{2(k-\beta)}{t+1}+\alpha+(s-l)(\alpha+\beta)\\
&=&n-\frac{2k}{t+1}+\beta(\frac{2}{t+1}-1)+[s-(l-1)](\alpha+\beta).
\end{eqnarray*}
Also note that $t\geq 2$ and $s\leq l-1$, we have $$j-i\leq n-\frac{2k}{t+1}=n-\frac{2[p(t+1)+d]}{t+1}=n-2p-\frac{2d}{t+1}.$$
The result follows.
\end{pf}

\vspace{3mm} As a consequence of Theorem \ref{Thm3} and Theorem \ref{Thm8}, we derive a formula for  projective dimension and regularity of the path ideals of a cycle, which generalize  the formula   obtained in \cite[Corollary 5.5]{AF1}

\begin{corollary}\label{cor2} Let $l,m,n,k,t $ be  integers as in Theorem \ref{Thm3} and Proposition \ref{prop2}. Suppose that
$\Delta_{m,l}(C_{n})$ $=\langle F_{1},\dots, F_{k}\rangle$ is the  path complex of the cycle  $C_{n}$ with standard labeling and $I_{m,l}(C_{n})=\mathcal{I}(\Delta_{m,l}(C_{n}))$ is the facet ideal of $\Delta_{m,l}(C_{n})$.
We can write $k$ as  $k=p(t+1)+d$, where $p\geq 0$ and $0\leq d\leq t$. Then

(1) The projective dimension of $R/I_{m,l}(C_{n})$   is given by
$$ pd\,(R/I_{m,l}(C_{n}))=\left\{\begin{array}{ll}
2p,&\mbox{if}\ d=0\\
2p+1,&\mbox{if}\ d\neq 0\\
\end{array}\right.$$

(2) The regularity of $R/I_{m,l}(C_{n})$   is given by
$$ reg\,(R/I_{m,l}(C_{n}))=\left\{\begin{array}{ll}
n-2p,&\mbox{if}\ d=0\\
n-2p-1,&\mbox{if}\ d\neq 0\\
\end{array}\right.$$

\end{corollary}
\begin{pf} If $t=1$, then   assertion is obvious from Theorem \ref{Thm3}, so we prove assertion is true for $t\geq 2$.

(1)  This follows from Theorem \ref{Thm3} and Theorem \ref{Thm5};

(2) By definition, $reg\,(R/I_{m,l}(C_{n}))=\mbox{max}\{j-i\ |\ \beta_{i,j}\neq 0\}$.
By Theorem \ref{Thm5} and Theorem \ref{Thm8} and the observation above, we have
$$ reg\,(R/I_{m,l}(C_{n}))=\left\{\begin{array}{ll}
n-2p,&\mbox{if}\ d=0\\
n-2p-1,&\mbox{if}\ d\neq 0.\\
\end{array}\right.$$
\end{pf}

\begin{corollary}\label{cor3} Let $l,m,n,k,t $ be  integers as in Proposition \ref{prop2}. Suppose that
$\Delta_{m,l}(C_{n})$ $=\langle F_{1},\dots, F_{k}\rangle$ is the  path complex of the cycle  $C_{n}$ with standard labeling and $I_{m,l}(C_{n})=\mathcal{I}(\Delta_{m,l}(C_{n}))$ is the facet ideal of $\Delta_{m,l}(C_{n})$.
Then $$depth\,(R/I_{m,l}(C_{n}))=reg\,(R/I_{m,l}(C_{n})).$$
\end{corollary}
\begin{pf}
 By Auslander-Buchsbaum formula, we get that
 $$depth\,(R/I_{m,l}(C_{n}))=n-pd\,(R/I_{m,l}(C_{n})).$$
It follows from Theorem \ref{Thm8}.
\end{pf}

 \vspace{20mm}
{\bf Acknowledgments}

 \hspace{3mm} I did  this work during my stay at Department of
Mathematics of University Duisburg-Essen, Germany.  I  would like to thank Professor J\"urgen Herzog for effective discussions.
I spent a memorable time at Essen, so I would like to express my hearty thanks  to Maja for hospitality.
I also wish to thank the hospitality of Department of
Mathematics of University Duisburg-Essen, Germany.

This research was  supported by the National Natural Science Foundation of
China (11271275) and by  Foundation of Jiangsu Overseas Research \& Training Program for
University Prominent Young \& Middle-aged Teachers and Presidents and
 by Foundation of the Priority Academic Program Development of Jiangsu Higher Education Institutions.

 \vspace{10mm}

Guangjun Zhu, School of Mathematic Science,Soochow University, Suzhou 215006, China

E-mail: zhuguangjun@suda.edu.cn

\end{document}